\documentclass[11pt]{article}
\usepackage{amsmath}
\usepackage{amssymb}
\usepackage{amsfonts}

\makeatletter

\@addtoreset{equation}{section}
\makeatother
\def\<{\langle}
\def\>{\rangle}

\def\div{{\rm div}\ }

\newtheorem{lem}{Lemma}[section]
\newtheorem{theo}{Theorem}[section]
\newtheorem{rem}{Remark}[section]
\newtheorem{pro}{Proposition}[section]

\makeatletter
   
   \@addtoreset{equation}{section}
\makeatother

\begin{document}
\title{\bf On the system of $2$-D elastic waves \\with critical space dependent damping}
\author{Ruy Coimbra Char\~ao\thanks{ruy.charao@ufsc.br}  \\{\small Department of Mathematics}, {\small Federal University of Santa Catarina} \\ {\small 88040-900, Florianopolis, SC,  Brazil} \\and\\ Ryo Ikehata\thanks{Corresponding author: ikehatar@hiroshima-u.ac.jp}\\{\small Department of Mathematics, Division of Educational Sciences}\\{\small Graduate School of Humanities and Soicial Sciences}\\{\small Hiroshima University}\\ {\small Higashi-Hiroshima 739-8524, Japan}}
\date{}
\maketitle
\begin{abstract}
We consider the system of elastic waves with critical space dependent damping $V(x)$.  We study the Cauchy problem for this model in the $2$-dimensional Euclidean space ${\bf R}^{2}$, and we obtain faster decay rates of the total energy as time goes to infinity. In the $2$-D case we do not have any suitable Hardy type inequality, so generally one has no idea to establish optimal energy decay. We develope a special type of multiplier method combined with some estimates brought by the $2$-D Newton potential belonging to the usual Laplacian $-\Delta$, not the operator $-a^2\Delta - (b^2-a^{2})\nabla \div$ itself. The property of finite speed propagation is important to get results for this system.
\end{abstract}

\section{Introduction}
\footnote[0]{Keywords and Phrases: Elastic waves; critical damping; $2$-D whole space; optimal decay estimates; multiplier method; Newton potential.}
\footnote[0]{2020 Mathematics Subject Classification. Primary 35L52, 35L05; Secondary 35B40, 35B45, 35J05.}

We consider the following dissipative system of elastic waves under effects of a Stokes damping term with variable spatial coefficient as follows:
\begin{align}
& u_{tt} - a^2\Delta u - (b^2-a^{2})\nabla\div u + V(x) u_t  = 0,\ \ \ (t,x)\in (0,\infty)\times {\bf R}^{2},\label{eqn}\\
& u(0,x)= u_0(x), \quad  u_{t}(0,x)= u_{1}(x),\ \ \ x\in{\bf R}^{2} ,\label{initial}
\end{align}
where $V(x) > 0$ is a damped bounded coefficient specified later, and the Lam\'e coefficients $a > 0$ and $b > 0$ satisfy $0 < a^{2} < b^{2}$.
The vector valued initial data shall be chosen as $u_{j} \in (C_{0}^{\infty}({\bf R}^{2}))^{2}$ ($j = 0,1$), and so a unique existence of smooth (classical) vector valued solution $u(t,x) \in {\bf R}^{2}$ to problem \eqref{eqn}-\eqref{initial} can be guaranteed by the standard Semigroup theory combined with the regularity argument (see Ikawa \cite{Ikawa}).
Furthermore, one has the energy identity
\begin{equation}\label{Energy-ident}
E_{u}(t) + \int_{0}^{t}\int_{{\bf R}^{2}}V(x)\vert u_{s}(s,x)\vert^{2}dxds \textbf{= }E_{u}(0),
\end{equation}
where
\[E_{u}(t) := \frac{1}{2}\int_{{\bf R}^{2}}\left(\vert u_{t}(t,x)\vert^{2} + a^{2}\vert \nabla u(t,x)\vert^{2} + (b^{2}-a^{2})(\div u(t,x))^{2}\right)dx\] 
represents the total energy with respect to the equation \eqref{eqn}.
\vspace{0.2cm}
\\

\underline{{\small Notation.}}\,For the vector $u = (u_{1}, u_{2}) \in {\bf R}^{2}$ ($u_{j} \in {\bf R}$, ($j = 1,2$)), we set $\vert u\vert := (\vert u_{1}\vert^{2} + \vert u_{2}\vert^{2})^{1/2}$. We set $\vert\nabla u\vert^{2} := \vert\nabla u_{1}\vert^{2} + \vert\nabla u_{2}\vert^{2}$ for the function $u(x) := (u_{1}(x),u_{2}(x))^{T}$, where $(p,q)^{T}$ represents a transposed vector of $(p,q)$. $\Delta u := (\Delta u^{(1)},\Delta u^{(2)})^{T}$, $u_{t} := (u_{t}^{(1)},u_{t}^{(2)})^{T}$, and $\nabla u\cdot \nabla v := \nabla u^{(1)}\cdot \nabla v^{(1)} + \nabla u^{(2)}\cdot \nabla v^{(2)}$ for $u:= (u^{(1)}, u^{(2)})^{T}$ and $v:= (v^{(1)}, v^{(2)})^{T}$. $\Vert w\Vert_{L^{p}}$ means the usual $L^{p}$-norm of $w \in L^{p}({\bf R}^{2})$ for $1 \leq p \leq \infty$.\\

On the real valued potential $V(x)$, one assumes that $V(x)$ is bounded, sufficently smooth in ${\bf R}^{2}$, and there exist a positive constants $V_0 > 0$ such that 
\begin{equation}\label{hypo-V}
0<\frac{V_0}{1+\vert x\vert} \leq V(x) \leq \Vert V\Vert_{L^{\infty}} < +\infty, \;\; x\in {\bf R}^2.
\end{equation}
That is, one considers the critical damping coefficient case. Additionally, for a technical reason we impose a rather stronger hypothesis on the initial data:
\begin{equation}\label{hypo-L}
{\rm supp}\, u_{0} \cup {\rm supp}\, u_{1} \subset \{\vert x\vert \leq L\}
\end{equation}
for some constant $L > 0$. By this assumption the corresponding solution $u(t,x)$ to problem \eqref{eqn}-\eqref{initial} satisfies the finite propagation speed property:
\begin{equation}\label{hypo-LT}
u(t,x) = 0\quad \vert x\vert > bt + L.
\end{equation} 
Incidentally, the coefficient $b$ corresponds to the speed of P-wave, and the coefficient $a$ is the speed of S-wave.

\vspace{0.2cm}
First, let's get some background on the issue.\\
For scalar-valued wave equations with $a=b$, that is, 
\[u_{tt} - a^2\Delta u + V(x) u_t  = 0,\]
when the friction coefficient is $V(x) \sim \vert x\vert^{-\alpha}$ ($\vert x\vert \to \infty$) is well studied, and in this case $\alpha = 1$ is said to be critical damping, $\alpha < 1$ is effective damping, and $\alpha > 1$ is non-effective friction, and the problems of identifying the best decay rate of total/local energy together with non-decay property have been well studied through their corresponding initial value problems in the whole space or exterior mixed problems. It should be noted, however, that in the case of two-dimensional space, the treatment of such problems is somewhat difficult to solve because of the inadequacy of good inequalities such as Hardy-type inequalities. These references include \cite{ITY, Soba} for $\alpha = 1$, \cite{Mo} for $\alpha > 1$, \cite{TY, W} for $\alpha < 1$, and the references therein. While, in \cite{Nakao-1, ike, D, LSM} the authors studied the energy decay problem in the case where the friction coefficient may be degenerate but truly floating in the spatial distance, that is, $V(x) \geq 0$ for all $x \in {\bf R}^{n}$, and $V(x) \geq V_{0} > 0$ for $x \in {\bf R}^{n}$ satisfying $\vert x\vert \geq L$ for large $L \gg 1$.

On the flip side, when considering the energy decay problem of the relevant problem \eqref{eqn}-\eqref{initial} of elastic waves, there does not seem to be much prior research on the case with a friction term $V(x)$ that depends only on spatial variables. The main ones include \cite{RR} and \cite{Ruy-Ike-2011}. For the scattering and/or local energy decay results to the damped and/or undamped elastic waves, one can refer the reader to \cite{RC, SS, K, Ka, DDK}. Incidentally, topics such as the behavior of elastic wave solutions that incorporate different types of friction effects (e.g., structural friction) and energy decay are also studied in detail in \cite{C, R, T, WCL} and the references therein. By the way, the best energy decay results for Eq. \eqref{eqn} in the two-dimensional exterior domain are published in \cite{H}, but it should be noted that since \cite{H} uses Hardy's inequality, it is an exterior domain-only result (with minor modifications), and the treatment of the case of the entire two-dimensional space is completely unresolved.

Here we aim to identify the precise rate of decay of the total energy of elastic waves that incorporate critical friction in two spatial dimensions with relatively few useful tools and, therefore, somewhat more difficult to handle. To this end, we further modify the useful method due to \cite{IM}, which modifies the Morawetz method \cite{mora}, to derive the essential inequalities at the core of the discussion. The important part is that it reduces to the estimate of the solution $h(x)$ of the two-dimensional (vector-valued) Poisson equation such as 
\[-\Delta h(x) = u_{1}(x) + V(x)u_{0}(x).\]

Our main result reads as follows. As a result, unlike the 2-D exterior domain problem with the Dirichlet null boundary condition such as \cite{H}, a $\log t$-growth correction is needed for the decay rate. This would be a phenomenon unique to the whole space unlike the exterior domain case.
\begin{theo}\label{theorem-1}
Assume that the bounded function $V \in C^{\infty}({\bf R}^{2})$ satisfies \eqref{hypo-V} and $u_{0} \in (C_{0}^{\infty}({\bf R}^{2}))^{2}$ and $u_{1} \in (C_{0}^{\infty}({\bf R}^{2}))^{2}$ satisfies \eqref{hypo-L} for some $L > 0$. Then, the smooth solution $u(t,x)$ to problem \eqref{eqn}-\eqref{initial} satisfies
\[E_{u}(t) = O(t^{-2}\log t)\quad (t \to \infty),\]
in the case when $V_{0} > 2b$,  
\[E_{u}(t) = O(t^{-\frac{V_{0}}{b}+\delta}\log t)\quad (t \to \infty),\]
for any small $\delta > 0$ in case of $b < V_{0} \leq 2b$.
\end{theo}
\begin{rem}{\rm Similar results would be obtained if the regularity conditions for $V(x)$ or initial values $u_{j}$ ($j = 0,1$) were relaxed further. What is important is the fact that new results are obtained even under such "typical" conditions as is assumed in the first stage. The generalization of regularity is left to the reader's interest. Of course, the condition of compact support for the initial value cannot be removed.}
\end{rem}

The proof of this theorem will be done in the next sections. As a byproduct of the proof of the main result (see Remark \ref{ike12-1}), we obtain the best estimate rating from the following top. This is also unprecedented in elastic waves, so it is noted.
\begin{theo}\label{theorem-2}
Assume that $V(x) = 0$, and $u_{0} \in (C_{0}^{\infty}({\bf R}^{2}))^{2}$ and $u_{1} \in (C_{0}^{\infty}({\bf R}^{2}))^{2}$ satisfies \eqref{hypo-L} for some $L > 0$. Then, the smooth solution $u(t,x)$ to problem \eqref{eqn}-\eqref{initial} with $V(x) = 0$ satisfies
\[\int_{{\bf R}^{2}}\vert u(t,x)\vert^{2}dx \leq C\left(\Vert u_{0}\Vert_{L^{2}}^{2} + \Vert u_{1}\Vert_{L^{\infty}}^{2} + \Vert u_{1}\Vert_{L^{1}}^{2}\log(2L+bt)\right),\quad (t \gg 1)\]
with some universal constant $C > 0$.
\end{theo}
\begin{rem}{\rm The above growth order appears to be the best in two dimensions when compared to the best growth estimate in logarithmic order for linear waves (i.e., $a = b$) obtained in \cite{JHDE-ike}. Of course, there is still the unsolved problem of obtaining an growth estimate of logarithmic order from below for elastic waves. By the way, for the explicit representations of solution for $3$-D elastic waves one can cite \cite{A, F}.}
\end{rem}


\section{Main Estimates}
To get the first fundamental estimates we use the multiplier $f(t)u_t + g(t)u$, where $u = u(t,x)$ is the vector valued solution of \eqref{eqn}--\eqref{initial}, and the smooth positive functions $f(t)$ and $g(t)$ will be specified later. Then one obtains
\begin{lem}
\begin{equation}\label{energy-2}
\int_{{\bf R}^2}\big[ u_{tt} - a^2\Delta u - (b^2-a^{2})\nabla \div u + V(x) u_t\big]\cdot\big[f(t)u_t + g(t)u\big]dx = \frac{d}{dt} e(t) + F(t) = 0,
\end{equation}
for $t \in (0,\infty)$, where $x\cdot y$ denotes the usual inner product between two vectors $x \in {\bf R}^{2}$ and $y\in {\bf R}^{2}$,    
\begin{equation}\label{def.E}
 e(t) = \int_{{\bf R}^2}\Big[\frac{f(t)}{2}\big[\; |u_t|^2 + a^2|\nabla u|^2 + (b^2-a^2)(\div u)^2 \big] + g(t)u\cdot u_t + \big[ V(x)g - g_t\big] \frac{|u|^2}{2} \Big]dx,   
 \end{equation}
 and
 \begin{align}\label{def.F}
F(t) =& \int_{{\bf R}^2}\Big[-\frac{f(t)}{2}\big[\; |u_t|^2 + a^2|\nabla u|^2 + (b^2-a^{2})(\div u)^2 \big] + \big[g_{tt}  - g_tV(x)  \big]\frac{|u|^2}{2}\nonumber\\
+& f(t)V(x)\frac{|u_t|^2}{2} + \Big[ g(t) \big[a^2|\nabla u|^2 + (b^2-a^{2})(\div u)^2 -|u_t|^2  \big] \Big]dx.  
 \end{align}
 \end{lem}
 
The proof of this lemma is simple and is similar to the same identity for the wave equation developed in \cite{ITY} and \cite[Lemma 2.1]{H} (see also \cite{Ruy-Ike-2011} for elastic waves).
 
Now, by rewriting $F(t)$, one can get
\begin{align}\label{ike3}
  F(t) = &\frac{1}{2} \int_{{\bf R}^2} \Big(\big[ 2V(x)f - f_t -2g\big]|u_t|^2
+ (2g-f_t) \big[a^2|\nabla u|^2 + (b^2-a^{2})(\div u)^2 \big] \Big)dx \nonumber\\
&+ \frac{1}{2}\int_{{\bf R}^2}(g_{tt}-V(x)g_t)|u|^2 dx.
\end{align}
 
Next we set
$$\Omega(t) := \{x \in {\bf R}^2: |x| \leq L+bt\}$$
defined for $t\geq 0$, and let us consider the solution in this region $\Omega(t)$ because of the property of finite speed of propagation \eqref{hypo-LT}.

The lemmas below give important inequalities.
\begin{lem}\label{lema2}
    Let $f(t)$ and $g(t)$ be smooth functions satisfying\\
    {\rm (i)}\,$2fV(x)-f_t-2g \geq 0, \;\; x\in \Omega(t)$, \\
    {\rm (ii)}\, $2g-f_t \geq 0$\\
for all $t \geq t_{0} \geq 0$ with some fixed $t_{0} \geq 0$. Let $u(t,x)$ be the smooth solution to problem \eqref{eqn}--\eqref{initial}. Then the following estimate holds for all $t\geq t_0 \geq 0$:
\begin{equation}\label{ike1}
\frac{d}{dt}\big[f(t)E_u(t)+g(t)\int_{{\bf R}^{2}}u\cdot u_t dx\big] \leq \frac{1}{2}\int_{{\bf R}^2}\big[V(x)g_t-g_{tt}\big]|u|^2 dx + \frac{d}{dt}\int_{{\bf R}^{2}}\big[ g_t-V(x)g\big]\frac{|u|^2}{2}dx,
\end{equation}
where $E_u(t)$ is given by \eqref{Energy-ident}.
\end{lem}
The proof of this lemma follows by using the assumptions (i) and (ii) of Lemma \ref{lema2}, \eqref{ike3} and identity \eqref{energy-2}.

\begin{lem}\label{lema3}
In addition to the conditions {\rm (i)} and {\rm (ii)} of Lemma {\rm \ref{lema2}} assumed for $f(t)$ and $g(t)$, let us suppose the following assumptions:\\
{\rm (iii)}\,$-g_{tt} \leq \frac{C_1}{1+t}$,\\
{\rm (iv)}\,$g_t(t) \leq C_2$,\\
{\rm (v)}\,$g_t-V(x)g \leq C_3,\;\;  x \in {\bf R}^2$,\\
for $t\geq t_0 \geq 0$, where \,$C_1, C_2$ and $C_3$ are positive constants. Then there exists a constant $C(t_{0}) > 0$ depending on $t_{0} > 0$ such that
\[f(t)E_u(t) +g(t)\int_{{\bf R}^{2}}(u(t,x)\cdot u_t(t,x))dx\]
\[\leq C(t_{0}) + \frac{C_3}{2}\int_{{\bf R}^2} |u(t,x)|^2 dx +\frac{C_2}{2}\int_{t_0}^t\int_{{\bf R}^2}V(x)|u(s,x)|^2dxds\]
\begin{equation}\label{ike2}
 +\frac{C_1}{2V_0}\int_{t_0}^t\int_{{\bf R}^2}\frac{V_0}{1+s}|u(s,x)|^2dxds \quad (t \geq t_{0}),
\end{equation}
where $V_{0} > 0$ is the constant defined in \eqref{hypo-V}. \\
{\it Proof.}\, Integrating both sides of \eqref{ike1} over $[t_{0}, t]$ and using the conditions {\rm (iii)}, {\rm (iv)} and {\rm (v)}, the unique solution $u = u(t,x)$ of \eqref{eqn}--\eqref{initial} satisfies the following inequalities:
\end{lem}
\[f(t)E_u(t) +g(t)\int_{{\bf R}^{2}}(u(t,x)\cdot u_t(t,x))dx\]
\[\leq f(t_0)E_u(t_0)+g(t_0)\int_{{\bf R}^{2}}\left( u(t_0,x)\cdot u_t(t_0,x)\right)dx + \frac{C_3}{2}\int_{{\bf R}^2} |u(t,x)|^2 dx\]
\[-\frac{1}{2}\int_{{\bf R}^2} \big[ g_t(t_0) -V(x)g(t_0) \big] |u(t_0,x)|dx+\frac{C_2}{2}\int_{t_0}^t\int_{{\bf R}^2}V(x)|u(s,x)|^2dxds\]
\[+\frac{C_1}{2}\int_{t_0}^t\int_{{\bf R}^2}\frac{1}{1+s}|u(s,x)|^2dxds \]
\[\leq f(t_0)E_u(t_0)+g(t_0)\int_{{\bf R}^{2}}\left\vert\big( u(t_0,x)\cdot u_t(t_0,x)\big)\right\vert dx+\frac{1}{2}\int_{{\bf R}^2}\left\vert \big[ g_t(t_0) -V(x)g(t_0) \big]\right\vert |u(t_0,x)|dx\]
\[+\frac{C_3}{2}\int_{{\bf R}^2} |u(t,x)|^2 dx+\frac{C_2}{2}\int_{t_0}^t\int_{{\bf R}^2}V(x)|u(s,x)|^2dxds\]
\[+\frac{C_1}{2V_0}\int_{t_0}^t\int_{{\bf R}^2}\frac{V_0}{1+s}|u(s,x)|^2dxds\]
with $V_0 >0$ defined in \eqref{hypo-V} . If we set
\[C(t_{0}) := f(t_0)E_u(t_0)+g(t_0)\int_{{\bf R}^{2}}\left\vert\big( u(t_0,x)\cdot u_t(t_0,x)\big)\right\vert dx\]
\begin{equation}\label{ike}
+\frac{1}{2}\int_{{\bf R}^2}\left\vert\big[ g_t(t_0) -V(x)g(t_0) \big]\right\vert |u(t_0,x)|dx,
\end{equation}
the desired estimate can be derived.
\hfill
$\Box$

Let us choose two functions $f(t)$ and $g(t)$ such as
\begin{equation}\label{f-g-1}
f(t) = (1+t)^2, \quad g(t) = 1+t
\end{equation}
for the case $V_0>2b$, and in the case  $b<V_0 \leq 2b$
\begin{equation}\label{f-g-2}
f(t) = (1+t)^{\frac{V_0}{b}-\delta}, \quad g(t) = \frac{V_{0}-b\delta}{2b}(1+t)^{\frac{V_{0}}{b}-1-\delta}
\end{equation}
with arbitrarily fixed $\delta>0$.

 Now we are in a position to prepare the next lemma (see \cite[Lemmas 2.2 and 2.3]{H}).

 \begin{lem}\label{i-v}
Let $f(t)$ and $g(t)$ defined by in \eqref{f-g-1} or \eqref{f-g-2}. Then, the conditions {\rm (i})--{\rm (v)} in Lemmas {\rm \ref{lema2}} and {\rm \ref{lema3}} hold on the region $\Omega(t)$ for $t\geq t_0 \gg 1$.
\end{lem}

The proof of this lemma is similar to the same one developed in previous works \cite{ITY} and /or \cite{H}, so we shall omit its proof. Incidentally, in order to prove Lemma \ref{i-v} one must take $t_{0} > 0$ sufficiently large as was discussed in \cite{H}, there the condition \eqref{hypo-V} was also used.

The consequence of Lemma \ref{i-v} and estimate \eqref{ike2}, for large $t_0>0$ is that
\begin{align}\label{ike4}
f(t)E_u(t) & +g(t)\int_{{\bf R}^{2}}(u(t,x)\cdot u_t(t,x))dx \leq  C(t_{0}) 
+\frac{C_2}{2}\int_{t_0}^t\int_{{\bf R}^2}V(x)|u(s,x)|^2dxds \\
&+\frac{C_3}{2}\int_{{\bf R}^2} |u(t,x)|^2 dx
+\frac{C_1}{2V_0}\int_{t_0}^t\int_{{\bf R}^2}\frac{V_0}{1+s}|u(s,x)|^2dxds. \nonumber
\end{align}


In the sequel we need to estimate the last integral in the above estimate. For the $2$-D case, since the Hardy type inequality never holds, we need to use a new idea to obtain such an estimate. This part is a main contribution throughout this paper. For that, we rely on the estimate for the $2$-D Newton potential. We use the idea due to \cite{ITY}, at this moment, to estimate the quantities indicated in the right hand side of \eqref{ike4}. For this, we set (see \cite{IM} for its origin of this idea)
\begin{equation}\label{ike7}
v(t,x) := \int_{0}^{t}u(s,x)ds,
\end{equation} 
for the solution $u(t,x)$. Then, in terms of the new function $v(t,x)$ we find that
\begin{align}
& v_{tt} - a^2\Delta v - (b^2-a^{2})\nabla\div v + V(x) v_t  = u_{1}(x) + V(x)u_{0}(x),\ \ \ (t,x)\in (0,\infty)\times {\bf R}^{2},\label{veqn}\\
& v(0,x)= 0, \quad  v_{t}(0,x)= u_{0}(x),\ \ \ x\in{\bf R}^{2}.\label{vinitial}
\end{align}
By developing the energy estimate to problem \eqref{veqn}-\eqref{vinitial} one has the energy identity:
\[\frac{1}{2}\int_{{\bf R}^{2}}\vert v_{t}(t,x)\vert^{2}dx + \frac{a^{2}}{2}\int_{{\bf R}^{2}}\vert \nabla v(t,x)\vert^{2}dx + \frac{(b^{2}-a^{2})}{2}\int_{{\bf R}^{2}}(\div v(t,x))^{2}dx \]
\begin{equation}\label{ike5}
+ \int_{0}^{t}\int_{{\bf R}^{2}}V(x)\vert v_{s}(s,x)\vert^{2}dxds = \frac{1}{2}\int_{{\bf R}^{2}}\vert u_{0}(x)\vert^{2}dx + \int_{{\bf R}^{2}}(u_{1}(x)+V(x)u_{0}(x))\cdot v(t,x)dx.
\end{equation}

In order to estimate the last term of \eqref{ike5} we define a vector valued function $h(x)$ by
\begin{equation}\label{ike6}
h(x) := -\frac{1}{2\pi}\int_{{\bf R}^{2}}(\log\vert x-y\vert)(u_{1}(y)+V(y)u_{0}(y))dy.
\end{equation} 
Incidentally, if we write $h(x) = (h^{1}(x),h^{2}(x))^{T}$, and $u_{j}(x) = (u^{1}_{j}(x), u^{2}_{j}(x))^{T}$ for $x \in {\bf R}^{2}$ and ($j = 0,1$), \eqref{ike6} means
\[h^{k}(x) := -\frac{1}{2\pi}\int_{{\bf R}^{2}}(\log\vert x-y\vert)(u_{1}^{k}(y)+V(y)u_{0}^{k}(y))dy\quad (k = 1,2).\]
It is known (see e.g., Evans \cite[page 23, Theorem 1]{E}) from the regularity on the initial data that $h \in C^{2}({\bf R}^{2},{\bf R}^{2})$ and the function $h(x)$ satisfies
\begin{equation}\label{ike6-1}
-\Delta h(x) = u_{1}(x)+V(x)u_{0}(x)\quad (x \in {\bf R}^{2})
\end{equation}
in the classical sense. Then, we get the following lemma.

\begin{lem}\label{N3}
The function $h(x)$ defined in \eqref{ike6} satisfies
\[\vert x\vert \vert \nabla h(x)\vert \leq C\Vert u_{1} + V(\cdot)u_{0}\Vert_{L^{1}}\]
for $\vert x \vert \geq 2L$, where $C > 0$ is a constant.
\end{lem}
{\it Proof of Lemma \ref{N3}.}\,\,Since it is sufficient to check each component individually of the vector-valued solution $h(x)$ and the initial data, it is enough to check $h(x)$ as a real-valued function together with initial data from the beginning. 

First of all, note that under the assumption $\vert x\vert \geq 2L$ we can get
\[\vert x-y\vert \geq \vert x\vert - L \geq L,\quad \vert x\vert -L \geq \frac{1}{2}\vert x\vert\]
for $y \in {\bf R}^{2}$ satisfying $\vert y\vert \leq L$. After simple elementary computations on \eqref{ike6} and the compact support assumption on the initial data $u_{j}$ ($j = 0,1$) one can get the estimate
\[\vert\nabla h(x)\vert \leq \frac{1}{2\pi}\int_{{\bf R}^{2}}\frac{\vert u_{1}(y) + V(y)u_{0}(y)\vert}{\vert x-y\vert}dy\]
\[= \frac{1}{2\pi}\int_{\vert y\vert \leq L}\frac{\vert u_{1}(y) + V(y)u_{0}(y)\vert}{\vert x-y\vert}dy \leq \frac{1}{2\pi}\frac{1}{\vert x\vert-L}\int_{{\bf R}^{2}}\vert u_{1}(y) + V(y)u_{0}(y)\vert dy\]
\[\leq \frac{1}{\pi}\frac{1}{\vert x\vert}\Vert u_{1}+V(\cdot)u_{0}\Vert_{L^{1}},\quad \vert x\vert \geq 2L,\]
which implies the desired estimate with a constant $C := \frac{1}{\pi}$. Final desired estimates for the vector-valued $h(x)$ can be easily derived with some universal constant $C > 0$.   

\hfill
$\Box$


\begin{lem}\label{ike8}\,For any $\varepsilon > 0$ there exists a constant $C_{\varepsilon}>0$ such that
\[\left\vert\int_{{\bf R}^{2}}(u_{1}(x)+V(x)u_{0}(x))\cdot v(t,x) dx\right\vert \leq C_{\varepsilon}\int_{\vert x\vert \leq 2L+bt}\vert\nabla h(x)\vert^{2}dx + \varepsilon\int_{{\bf R}^{2}}\vert\nabla v(t,x)\vert^{2}dx \quad (t \geq 0).\]
\end{lem}
\begin{rem}{\rm If we use the original idea due to Morawetz \cite{mora}, we must solve more complicated equation 
\[-a^2\Delta h(x) - (b^2-a^{2})\nabla\div h(x) = -(u_{1}(x)+V(x)u_{0}(x)),\]
and must treat the following function $v(t,x)$ with its auxiliary function $h(x)$ in place of \eqref{ike7}:
\[v(t,x) := h(x) + \int_{0}^{t}u(s,x)ds.\]
However, in our new idea it suffices to solve the simple Poisson equation
\[-\Delta h(x) = u_{1}(x)+V(x)u_{0}(x),\]
and we just treat \eqref{ike7}.}
\end{rem}
{\it Proof of Lemma \ref{ike8}}\, It follows from the finite speed of propagation property for $v(t,x)$, \eqref{ike6-1} and the integration by parts that
\[\left\vert \int_{{\bf R}^{2}}(u_{1}(x) + V(x)u_{0}(x))\cdot v(t,x)dx\right\vert = \left\vert \int_{\vert x\vert \leq 2L + bt}(u_{1}(x)+ V(x)u_{0}(x))\cdot v(t,x)dx \right\vert\]
\[= \left\vert -\int_{\vert x\vert \leq 2L + bt}\Delta h(x)\cdot v(t,x)dx \right\vert = \left\vert \int_{\vert x\vert \leq 2L + bt}\nabla h(x)\cdot\nabla v(t,x)dx \right\vert\]
\[\leq 2\int_{\vert x\vert \leq 2L + bt}\vert\nabla h(x)\vert\vert\nabla v(t,x)\vert dx\]
\[\leq C_{\varepsilon}\int_{\vert x\vert \leq 2L + bt}\vert\nabla h(x)\vert^{2}dx + \varepsilon\int_{\vert x\vert \leq 2L + bt}\vert\nabla v(t,x)\vert^{2}dx\]
\begin{equation}\label{N5}
\leq C_{\varepsilon}\int_{\vert x\vert \leq 2L + bt}\vert\nabla h(x)\vert^{2}dx + \varepsilon\int_{{\bf R}^{2}}\vert\nabla v(t,x)\vert^{2}dx
\end{equation} 
with some parameter $\varepsilon > 0$ and a constant $C_{\varepsilon} > 0$ depending only on $\varepsilon > 0$, where one has just used the facts that $v(t,x) = 0$, and $\vert\nabla h(x)\vert < +\infty$ for $\vert x\vert = 2L + bt$ for the boundary integral. Note that 
\[\int_{\vert x\vert \leq 2L + bt}\vert\nabla h(x)\vert^{2}dx < +\infty\]
for each $t \geq 0$ because of $h \in C^{2}({\bf R}^{2}, {\bf R}^{2})$. These arguments imply the desired estimate.

\hfill
$\Box$

\vspace{0.2cm}
Let us estimate the first term of the right hand side of Lemma \ref{ike8} by using Lemma \ref{N3}. For this purpose we set
\begin{equation}\label{ruy-1}
I_{h} := \int_{\vert x\vert \leq 2L}\vert\nabla h(x)\vert^{2}dx.
\end{equation}
First, we see that
\begin{equation}\label{N8}
\int_{\vert x\vert \leq 2L + bt}\vert\nabla h(x)\vert^{2}dx = I_{h} + \int_{2L \leq \vert x\vert \leq 2L + bt}\vert\nabla h(x)\vert^{2}dx.
\end{equation}
Thus, by using polar coordinates from Lemma \ref{N3} it follows that 
\[\int_{2L \leq \vert x\vert \leq 2L + bt}\vert\nabla h(x)\vert^{2}dx \leq C^{2}\Vert u_{1}+V(\cdot)u_{0}\Vert_{L^{1}}^{2}\int_{2L \leq \vert x\vert \leq 2L + bt}\frac{1}{\vert x\vert^{2}}dx\]
\[= 2\pi C^{2}\Vert u_{1}+ V(\cdot)u_{0}\Vert_{L^{1}}^{2}\int_{2L}^{2L + bt}\frac{r}{r^{2}}dr\]
\[= 2\pi C^{2}\Vert u_{1}+ V(\cdot)u_{0}\Vert_{L^{1}}^{2}\left(\log(2L+bt)-\log(2L)\right)\]
\begin{equation}\label{N9}
\leq 2\pi C^{2}\Vert u_{1}+ V(\cdot)u_{0}\Vert_{L^{1}}^{2}\log(2L+bt)\quad (t \geq 0). 
\end{equation}
Therefore, by \eqref{N8} and \eqref{N9} one has arrived at the crucial estimate:
\begin{lem}\label{N10}
It holds that
\[\int_{\vert x\vert \leq  2L + bt}\vert\nabla h(x)\vert^{2} dx \leq I_{h} + 2\pi C^{2}\Vert u_{1}+ V(\cdot)u_{0}\Vert_{L^{1}}^{2} \log (2L+bt) \quad (t \geq 0). \]
\end{lem}

Finally, let us estimate $I_{h}$ in terms of the initial velocity $u_{1}+ V(\cdot)u_{0}$. Note that generally it holds that
\[\int_{{\bf R}^{2}}\vert\nabla h(x)\vert^{2}dx = +\infty.\] 
It is important to point out that the integral region of the quantity $I_{h}$ is localized to $\vert x\vert \leq 2L$ as is defined in \eqref{ruy-1}. Of course, we know that its integral value of $I_{h}$ is finite because of $h \in C^{2}({\bf R}^{2},{\bf R}^{2})$, but we will try to suppress it by a certain quantity of initial velocity from above.\\

Thus, it follows from \eqref{ike5} and Lemmas \ref{ike8} and \ref{N10} that
\[\frac{1}{2}\int_{{\bf R}^{2}}\vert v_{t}(t,x)\vert^{2}dx + \left(\frac{a^{2}}{2}-\varepsilon\right)\int_{{\bf R}^{2}}\vert \nabla v(t,x)\vert^{2}dx\]
\[+ \frac{(b^{2}-a^{2})}{2}\int_{{\bf R}^{2}}(\div v(t,x))^{2}dx + \int_{0}^{t}\int_{{\bf R}^{2}}V(x)\vert v_{s}(s,x)\vert^{2}dxds\]
\begin{equation}\label{ike9}
\leq \frac{1}{2}\int_{{\bf R}^{2}}\vert u_{0}(x)\vert^{2}dx + C_{\varepsilon}I_{h} + 2\pi C_{\varepsilon} C^{2}\Vert u_{1}+ V(\cdot)u_{0}\Vert_{L^{1}}^{2} \log (2L+bt).
\end{equation}
 
For the quantity $I_{h}$ one can get the following result.
\begin{lem}\label{ike11}\,It hods that
\[I_{h} = \int_{\vert x\vert \leq 2L}\vert\nabla h(x)\vert^{2}dx \leq C_{L}\Vert u_{1}+V(\cdot)u_{0}\Vert_{L^{\infty}}^{2},\quad (t \geq 0),\]
where $C_{L} > 0$ is a constant depending on $L > 0$.
\end{lem}
{\it Proof.}\,The check here is also initially performed when $h(x)$ and the initial value $u_{j}$ ($j = 0,1$) is a real valued, as follows. 

For an arbitrarily fixed $x_{0} \in \overline{B_{2L}(0)}$, from the definition of $h(x)$ in \eqref{ike6} one can easily arrive at the intermediate estimate:
\[\vert\nabla h(x_{0})\vert \leq \frac{1}{2\pi}\int_{{\bf R}^{2}}\frac{\vert u_{1}(y)+V(y)u_{0}(y)\vert}{\vert x_{0}-y\vert}dy = \frac{1}{2\pi}\int_{\vert y\vert \leq L}\frac{\vert u_{1}(y)+V(y)u_{0}(y)\vert}{\vert x_{0}-y\vert}dy \]
\[= \frac{1}{2\pi}\int_{\overline{B_{L}(x_{0})}}\frac{\vert u_{1}(x_{0}-z)+V(x_{0}-z)u_{0}(x_{0}-z)\vert}{\vert z\vert}dz\]
\begin{equation}\label{N11}
\leq \frac{\Vert u_{1}+V(\cdot)u_{0}\Vert_{L^{\infty}}}{2\pi}\int_{\overline{B_{L}(x_{0})}}\frac{1}{\vert z\vert}dz.
\end{equation}
It should be noticed that in the case of $x_{0} \in \overline{B_{2L}(0)}$, we have $\overline{B_{L}(x_{0})} \subset \overline{B_{4L}(0)}$ sufficiently. Thus, we get
\begin{equation}\label{N12}
\int_{\overline{B_{L}(x_{0})}}\frac{1}{\vert z\vert}dz \leq \int_{\overline{B_{4L}(0)}}\frac{1}{\vert z\vert}dz = 8\pi L.
\end{equation}
Put $x_{0}$ back into $x$ to get the following by \eqref{N11} and \eqref{N12}. 
\[\vert\nabla h(x)\vert \leq 4L\Vert u_{1}+V(\cdot)u_{0}\Vert_{L^{\infty}}\quad \forall x \in \overline{B_{2L}(0)}.\]
Therefore, one can obtain the desired estimate with some constant $C_{L} > 0$ depending only on $L$. By a simple calculation, the exact same estimate holds for vector-valued $h(x)$ and initial values $u_{j}$ ($j = 0,1$).
\hfill
$\Box$

By combining these estimates developed in \eqref{ike9} and Lemma \ref{ike11}, and by choosing $\varepsilon > 0$ small enough in \eqref{ike9}, one can arrive at the crucial estimate because of $v_{t} = u$. Note that at this stage one must assume $a > 0$ necessarily. 
\begin{pro}\label{ike12}\,There exists a constant $C > 0$ depending on $L > 0$ such that 
\[\int_{{\bf R}^{2}}\vert u(t,x)\vert^{2}dx + \int_{0}^{t}\int_{{\bf R}^{2}}V(x)\vert u(s,x)\vert^{2}dxds\]
\[\leq C\left(\Vert u_{0}\Vert_{L^{2}}^{2} + \Vert u_{1}+V(\cdot)u_{0}\Vert_{L^{\infty}}^{2} + \Vert u_{1}+V(\cdot)u_{0}\Vert_{L^{1}}^{2}\log(2L+bt)\right)\]
for all $t \geq 0$.
\end{pro}
\begin{rem}\label{ike12-1}{\rm By the results from \cite{JHDE-ike} for scalar valued wave equation with $a=b$ and $V(x) = 0$, it seems that the $\log$-order growth rate obtained in Proposition \ref{ike12} may be optimal in the two dimensional case. Note that $V(x) = 0$ case can be admitted at this stage, that is, for the solution $w(t,x)$ satisfying
\begin{align}
& w_{tt} - a^2\Delta w - (b^2-a^{2})\nabla\div w  = 0,\ \ \ (t,x)\in (0,\infty)\times {\bf R}^{2},\label{2eqn}\\
& w(0,x) = u_0(x), \quad  w_{t}(0,x) = u_{1}(x),\ \ \ x\in{\bf R}^{2} ,\label{2initial}
\end{align}
it follows from the proof of Proposition \ref{ike12} that
\[\int_{{\bf R}^{2}}\vert w(t,x)\vert^{2}dx \leq C\left(\Vert u_{0}\Vert_{L^{2}}^{2} + \Vert u_{1}\Vert_{L^{\infty}}^{2} + \Vert u_{1}\Vert_{L^{1}}^{2}\log(2L+bt)\right),\quad (t \gg 1)\]
with some constant $C > 0$. Incidentally, the discussion here is also directly related to the proof of Theorem \ref{theorem-2}.}
\end{rem}

Now, we use the assumption \eqref{hypo-V} (which means $V(x) > 0$). Indeed, it should be pointed out that from \eqref{hypo-V} and \eqref{hypo-LT} one has
\[\int_{t_{0}}^{t}\int_{{\bf R}^{2}}V(x)\vert u(s,x)\vert^{2}dxds = \int_{t_{0}}^{t}\int_{\vert x\vert \leq L+bs}V(x)\vert u(s,x)\vert^{2}dxds \]
\begin{equation}\label{ike13}
\geq V_{0}\int_{t_{0}}^{t}\int_{\vert x\vert \leq L+bs}(1+|x|)^{-1} \vert u(s,x)\vert^{2}dxds \geq C^{*}\int_{t_{0}}^{t}\int_{{\bf R}^{2}}\frac{V_{0}}{1+s}\vert u(s,x)\vert^{2}dxds,\quad t_{0} \gg 1,
\end{equation}
where 
\[C^{*} := \frac{1}{\max\{(1+L), b\}}.\]
Thus, the last term of \eqref{ike4} can be absorbed into the quantity $\displaystyle{\int_{t_{0}}^{t}\int_{{\bf R}^{2}}}V(x)\vert u(s,x)\vert^{2}dxds$, so that it follows from \eqref{ike4} and Proposition \ref{ike12} that
\[f(t)E_{u}(t) + g(t)\int_{{\bf R}^{2}}(u(t,x)\cdot u_{t}(t,x))dx \leq C(t_{0}) \]
\[+ \max\{\frac{C_{2}}{2},\frac{C_{3}}{2}\}C\left(\Vert u_{0}\Vert_{L^{2}}^{2} + \Vert u_{1}+V(\cdot)u_{0}\Vert_{L^{\infty}}^{2} + \Vert u_{1}+V(\cdot)u_{0}\Vert_{L^{1}}^{2}\log(2L+bt)\right)\]
\[+ \frac{C_{1}C}{2V_{0}^{2}C^{*}}\left(\Vert u_{0}\Vert_{L^{2}}^{2} + \Vert u_{1}+V(\cdot)u_{0}\Vert_{L^{\infty}}^{2} + \Vert u_{1}+V(\cdot)u_{0}\Vert_{L^{1}}^{2}\log(2L+bt)\right) \quad(t \geq t_{0}).\]
This implies the following lemma.
\begin{lem}\label{ike14} Undet the assumption as in Theorem {\rm \ref{theorem-1}} for $t \geq t_{0} \gg 1$ it holds that
\[f(t)E_{u}(t) + g(t)\int_{{\bf R}^{2}}(u(t,x)\cdot u_{t}(t,x))dx\]
\[\leq C(t_{0}) + C^{**}\left(\Vert u_{0}\Vert_{L^{2}}^{2} + \Vert u_{1}+V(\cdot)u_{0}\Vert_{L^{\infty}}^{2} + \Vert u_{1}+V(\cdot)u_{0}\Vert_{L^{1}}^{2}\log(2L+bt)\right)\]
with some universal constant $C^{**} > 0$ depending only on $V_{0}$, $L$, $a$ and $b$. 
\end{lem}

Let us finalize the proof of Theorem \ref{theorem-1}. \\
{\it Proof of Theorem \ref{theorem-1} completed.}\,First, by choosing $t_{0} > 0$ large enough, it follows from Proposition \ref{ike12} that
\[\int_{{\bf R}^{2}}\vert u(t,x)\vert^{2}dx \leq C\left(\Vert u_{0}\Vert_{L^{2}}^{2} + \Vert u_{1}+V(\cdot)u_{0}\Vert_{L^{\infty}}^{2} + \Vert u_{1}+V(\cdot)u_{0}\Vert_{L^{1}}^{2}\right)\log(2L+bt)\]
for all $t \geq t_{0}$. Set 
\[A_{0} := C\left(\Vert u_{0}\Vert_{L^{2}}^{2} + \Vert u_{1}+V(\cdot)u_{0}\Vert_{L^{\infty}}^{2} + \Vert u_{1}+V(\cdot)u_{0}\Vert_{L^{1}}^{2}\right).\]
Then, from Lemma \ref{ike14} and \eqref{Energy-ident} one can get the estimate
\[f(t)E_{u}(t) \leq 2g(t)\left(\int_{{\bf R}^{2}}\vert u(t,x)\vert^{2}dx\right)^{1/2}E_{u}(t)^{1/2} \]
\[+\,C(t_{0}) +\,C^{**}\left(\Vert u_{0}\Vert_{L^{2}}^{2} + \Vert u_{1}+V(\cdot)u_{0}\Vert_{L^{\infty}}^{2} + \Vert u_{1}+V(\cdot)u_{0}\Vert_{L^{1}}^{2}\log(2L+bt)\right)\]
\[\leq 2g(t)E_{u}(t)^{1/2}A_{0}^{1/2}\sqrt{\log(2L+bt)} +\,C(t_{0})\]
\[+\,C^{**}\left(\Vert u_{0}\Vert_{L^{2}}^{2} + \Vert u_{1}+V(\cdot)u_{0}\Vert_{L^{\infty}}^{2} + \Vert u_{1}+V(\cdot)u_{0}\Vert_{L^{1}}^{2}\log(2L+bt)\right)\]
\[\leq 2\frac{g(t)}{\sqrt{f(t)}}\sqrt{f(t)E_{u}(t)}A_{0}^{1/2}\sqrt{\log(2L+bt)} +\,C(t_{0})\]
\[+\,C^{**}\left(\Vert u_{0}\Vert_{L^{2}}^{2} + \Vert u_{1}+V(\cdot)u_{0}\Vert_{L^{\infty}}^{2} + \Vert u_{1}+V(\cdot)u_{0}\Vert_{L^{1}}^{2}\log(2L+bt)\right)\]
\[\leq 2A_{0}\frac{g(t)^{2}}{f(t)}\log(2L+bt) + \frac{1}{2}f(t)E_{u}(t) + \,C(t_{0})\]
\[+\,C^{**}\left(\Vert u_{0}\Vert_{L^{2}}^{2} + \Vert u_{1}+V(\cdot)u_{0}\Vert_{L^{\infty}}^{2} + \Vert u_{1}+V(\cdot)u_{0}\Vert_{L^{1}}^{2}\log(2L+bt)\right),\]
which implies the following estimate:
\[\frac{1}{2}f(t)E_{u}(t) \leq 2A_{0}\frac{g(t)^{2}}{f(t)}\log(2L+bt) + \,C(t_{0})\]
\[+\,C^{**}\left(\Vert u_{0}\Vert_{L^{2}}^{2} + \Vert u_{1}+V(\cdot)u_{0}\Vert_{L^{\infty}}^{2} + \Vert u_{1}+V(\cdot)u_{0}\Vert_{L^{1}}^{2}\log(2L+bt)\right),\]
so that one has the crucial estimate in terms of the functions $f(t)$ and $g(t)$:
\[E_{u}(t) \leq 4A_{0}\frac{g(t)^{2}}{f(t)^{2}}\log(2L+bt) + C^{**}\frac{2}{f(t)}\Vert u_{0}\Vert_{L^{2}}^{2} + \frac{2C(t_{0})}{f(t)}\]
\begin{equation}\label{ike15}
+\,C^{**}\frac{2}{f(t)}\left(\Vert u_{1}+V(\cdot)u_{0}\Vert_{L^{\infty}}^{2} + \Vert u_{1}+V(\cdot)u_{0}\Vert_{L^{1}}^{2}\log(2L+bt)\right),
\end{equation}
for $t \geq t_{0} \gg 1$.\\
Based on the estimate \eqref{ike15} one can get the desired first estimate
\[E_{u}(t) = O(t^{-2}\log t)\quad (t \to \infty),\]
in the case when $V_{0} > 2b$, by choosing $f(t) = (1+t)^{2}$ and $g(t) = 1+t$ (see \eqref{f-g-1}). Furthermore, in case of $b < V_{0} \leq 2b$, by taking $f(t) = (1+t)^{\frac{V_{0}}{b}-\delta}$ and $g(t) = \frac{V_{0}-b\delta}{2b}(1+t)^{\frac{V_{0}}{b}-1-\delta}$ (see \eqref{f-g-2}) one can get the "decay" estimate
\[E_{u}(t) = O(t^{-\frac{V_{0}}{b}+\delta}\log t)\quad (t \to \infty)\]
for small $\delta > 0$. These show the validity of our main theorem. Observe that the choice of functions $f(t)$ and $g(t)$ depends on the cases under the Lam\'e coefficient $b$ such that $V_{0} > 2b$ or $b \leq V_{0} \leq 2b$ because these conditions are necessary to the proof of Lemma \ref{i-v}.   
\hfill
$\Box$
\begin{rem}{\rm See \cite{Ruy-Ike-2011} for results when $0 < V_{0} \leq b$.}
\end{rem}


\vspace{0.5cm}
\noindent{\em Acknowledgement.}
\smallskip The work of the first author (R. C. Char\~ao) was partially supported by Print/Capes - Process 88881.310536/2018-00. 


\vspace{0.5cm}

{\bf Data availability}\\
\noindent
No data was used for the research described in the article.

{\bf Declarations}\\
\noindent
{\bf Conflict of interest} The authors have no conflict of interest, financial or otherwise.

\end{document}